\documentclass[12pt,reqno]{amsart}

\usepackage{amssymb}
\usepackage{amsmath}
\usepackage{latexsym}
\usepackage{amsthm}
\usepackage{psfig}

\theoremstyle{plain}
\newtheorem{theorem}{Theorem}[section]

\newtheorem{corollary}{Corollary}[section]
\newtheorem{lemma}{Lemma}[section]
\newtheorem{remark}{\it Remark}[section]
\theoremstyle{definition}
\newtheorem{definition}{Definition}[section]

\newcommand{\rot}{\mathop\mathrm{rot}}

\renewcommand{\div}{\mathop\mathrm{div}}
\newcommand{\vphi}{\varphi}
\newcommand{\ve}{\varepsilon}

\newcommand{\Zbb}{\ensuremath{\mathbb{Z}}}

\title
[  Sharp estimates for
the number of  degrees of freedom ]
 {Sharp estimates for
the number of  degrees of freedom for
 the  damped-driven 2D Navier--Stokes equations}
\author[A.A. Ilyin and E.S. Titi]{}

\begin{document}

\maketitle

\centerline{\scshape  Alexei A. Ilyin\footnote{\noindent Keldysh
Institute of
Applied Mathematics, Russian Academy of Sciences, Miusskaya Sq. 4, 125047
Moscow, Russia, E-mail: ilyin@spp.keldysh.ru } and
Edriss S. Titi\footnote
{Department of Mathematics and Department of Mechanical
and Aerospace Engineering, University of California,
 Irvine, California 92697, USA, E-mail: etiti@math.uci.edu. Also:
  Department of Computer Science and Applied Mathematics,
  Weizmann Institute of Science,
  P.O. Box 26,
  Rehovot, 76100, Israel,
  E-mail: edriss.titi@weizmann.ac.il 
}}

\medskip

\bigskip
\begin{quote}{\normalfont\fontsize{8}{10}\selectfont
{\bfseries Abstract.}
We derive upper bounds for the number of asymptotic degrees
(determining modes and  nodes)
of freedom for the  two-dimensional Navier--Stokes system and Navier-Stokes
system with damping.
In the first case we obtain the previously known estimates in
an explicit form, which are larger than the fractal dimension
of the global attractor. However, for the Navier--Stokes
system with damping
our estimates for the number of the determining modes and nodes
are comparable to the sharp estimates for the fractal dimension
of the global attractor. Our investigation of the damped-driven
2D Navier--Stokes system is inspired by the Stommel--Charney barotropic
model of ocean circulation where the damping represents the
Rayleigh friction. We remark that our results equally apply to
 the Stommel--Charney model.

\medskip


\medskip

\par}
\end{quote}

\setcounter{equation}{0}
\section{Introduction}\label{S:Intro}

In this paper we derive estimates for the number of determining
modes, nodes and other determining projections for the two-dimensional
Navier--Stokes system~(\ref{NSE}) and for the two-dimensional
damped-driven Navier--Stokes system~(\ref{NSE-damped}). The latter
system is inspired by the viscous Stommel--Charney barotropic
ocean circulation model~\cite{Char},\cite{P},\cite{Stom}:
\begin{equation}\label{Charney-Stommel}
\begin{aligned}
\partial_tu+\sum_{i=1}^2u^i\partial_iu+\mathbf{k}l\times u&=
-\mu u+\nu\Delta\,u-\nabla\,p +f,\\
\div u&=0,
\end{aligned}
\end{equation}
where the damping $\mu u$ represents the so-called
Rayleigh friction term in the ocean circulation model.
In recent years there has been some analytical study
of the Stommel--Charney model from the
dynamical systems point of view
(see, for instance,~\cite{BCT}, \cite{C-G-S-W},
\cite{Hauk}, \cite{Hauk-Titi},
\cite{I91}, \cite{I-M-T}, \cite{S}, \cite{W}).
For the sake of clarity in our  mathematical presentation
and in order to make a  straight forward comparison to the
physical literature about the 2D turbulence we focus here 
on the Navier--Stokes
 system~(\ref{NSE})
and the damped-driven Navier--Stokes
 system ~(\ref{NSE-damped}).
It worth stressing, however, that our results concerning
the system~(\ref{NSE-damped}), especially about determining modes,
equally applies to the barotropic model~(\ref{Charney-Stommel}),
which we report in a forthcoming paper.

Our paper is organized as follows.
In section~\ref{S:N-S} we obtain previously known estimates
reported in~\cite{Cockburn-Jones-Titi-1},
\cite{Cockburn-Jones-Titi} and \cite{Jones-Titi2}
for the number of determining modes, nodes
and other determining projections
for the two-dimensional space-periodic Navier--Stokes system
which are linear with respect to the Grashof number. 
We use, however, the scalar  vorticity formulation,
which makes it possible to give all the estimates
and constants in
an explicit form. Furthermore,
the dependence on the aspect ratio of the periodic domain
is explicitly singled out. For the Dirichlet boundary conditions we
obtain an explicit estimate for the number of determining modes
which is quadratic with respect to the Grashof number
(see, for example, \cite{H-Ti}, \cite{Ju}). It is worth
mentioning that these 
best known estimates for the number of explicit determining modes and 
nodes of the 2D Navier--Stokes equations without damping  are 
 still much larger than the dimension of the 
global attractor.  

In section~\ref{S:damped N-S} we consider the
Navier-Stokes system with damping, subject to
periodic boundary conditions and the stress-free boundary conditions,
and obtain
estimates for the number of the determining modes and nodes
that are  of the same order as of the sharp 
estimates for the fractal
dimension of the global attractor~\cite{I-M-T}.
These remarkable estimates are extensive, that is, depend linearly on the
area of the spatial domain, which is consistent with the
physical intuition. We remark, again, that such an observation is not
known to exist in the case of the 2D Navier--Stokes equations
without damping. 

Finally, in the Appendix in section~\ref{S:Appen}
we prove some auxiliary inequalities, namely,
we derive sharp constant in the Agmon inequality on the
two-dimensional torus and prove a variant of the embedding theorem
reported in~\cite{Jones-Titi2}.

Since we are interested in obtaining explicit bounds about the constants
involved in our asymptotic estimates for the numbers of degrees of freedom,
we focus in this paper on the notions of determining modes and nodes.
However, it is worth stressing that our asymptotic estimates, in terms
of the physical parameters, are valid for other determining
functionals and projections, as it has
been  demonstrated in \cite{Cockburn-Jones-Titi-1},
\cite{Cockburn-Jones-Titi},   (see also \cite{Chueshov}),
with constants that may vary depending on the underlying chosen
determining functionals.

\setcounter{equation}{0}
\section{Determining modes and nodes for two-dimensional
Navier--Stokes equations}\label{S:N-S}
\subsection*{Dirichlet boundary conditions}
We consider in this section the two-dimensional Navier--Stokes system
\begin{equation}\label{NSE}
\begin{aligned}
\partial_tu+\sum_{i=1}^2u^i\partial_iu&=
\nu\Delta\,u-\nabla\,p +f,\\
\div u&=0,\quad
u(0)=u_0,
\end{aligned}
\end{equation}
where $u$ is the velocity vector field satisfying Dirichlet
boundary conditions $u\vert_{\partial\Omega}=0$, $p$ is the pressure,
and $\nu>0$ is the kinematic viscosity. The right-hand side
$f=f(x,t)$ is given and the domain $\Omega$ is an
arbitrary open connected set in $\mathbb{R}^2$ with
finite measure $|\Omega|<\infty$.

We use the standard notation and facts from
the theory of Navier--Stokes equations
(see, for instance, \cite{CF88}, \cite{Lad}, \cite{TNS}, \cite{T})
and denote by $P$ the Helmholtz--Leray orthogonal projection in
$L_2(\Omega)^2$ onto the Hilbert space $H$
which is the closure in $L_2(\Omega)^2$ of the set of smooth
solenoidal vector functions with compact supports in~$\Omega$.
Applying $P$ to the first
equation in~(\ref{NSE}), we obtain
\begin{equation}\label{FNSE}
\partial_tu+B(u,u)+\nu Au=f, \qquad u(0)=u_0,
\end{equation}
where $A=-P\Delta$ is the Stokes operator,
$B(u,v)=P\bigl(\sum_{i=1}^2u^i\partial_iv\bigr)$ is the nonlinear term,
 and $f=Pf\in H$.

Next, we denote by
$\{\lambda_j\}_{j=1}^\infty$, $0<\lambda_1\le\lambda_2\le\dots$ and
$\{w_j\}_{j=1}^\infty$ the eigenvalues and the corresponding eigenfunctions
of the Stokes operator~$A$:
$A\,w_j=\lambda_jw_j\,$. The  asymptotic behavior
$\lambda_k\sim\frac{4\pi k}{|\Omega|}$ as $k\to\infty$
was established in~\cite{Metiv}, while in this work
we use the following explicit non-asymptotic lower bounds
 for the
eigenvalues~$\{\lambda_j\}_{j=1}^\infty$ (see~\cite{I96}):
\begin{equation}\label{lower-bounds}
\sum_{j=1}^m\lambda_j\ge\frac{\pi m^2}{|\Omega|}
\quad\text{and, consequently,}\quad\lambda_m\ge
\frac{\pi m}{|\Omega|}\,,\ m\ge1;\quad
\lambda_1\ge\frac{2\pi}{|\Omega|}\,.
\end{equation}

The Hilbert space $V=D(A^{1/2})$ is the space $H^1_0(\Omega)^2\cap H$
with norm
$$
\|u\|_{D(A^{1/2})}=\|\nabla\,u\|=\|\rot u\|,
$$
where $\|\cdot\|=\|\cdot\|_{L_2(\Omega)}$.
The nonlinear operator $B(v,v)$ satisfies the
well-known inequalities
(see, for instance,~\cite{CF88}, \cite{Lad}, \cite{TNS}, \cite{T})
\begin{equation}\label{est-for-B}
\aligned
&|(B(v,v),u)|=|(B(v,u),v)|\,\le\, c_1\|v\|\|\nabla\,v\|\|\nabla\,u\|,\\
&|(B(u,v),w)|\,\le\,
c_2
\|u\|^{1/2}\|\nabla\,u\|^{1/2}\|\nabla\,v\|\|w\|^{1/2}\|\nabla\,w\|^{1/2},
\endaligned
\end{equation}
where it was shown in~\cite{CI} that
\begin{equation}\label{est-for-cb}
    c_1=:c_b\le\left(\frac8{27\pi}\right)^{1/2}
    \qquad c_2\le\sqrt{2}c_b\le\left(\frac{16}{27\pi}\right)^{1/2}.
\end{equation}

Let $u$ and $v$ be the solutions of the Navier--Stokes
equations
\begin{equation}\label{two-NS-eq}
\begin{aligned}
&\partial_tu+B(u,u)+\nu Au=f, \qquad u(0)=u_0,\\
&\partial_tv+B(v,v)+\nu Av=g, \qquad v(0)=v_0,
\end{aligned}
\end{equation}
where $f,g\in L_\infty(0,\infty;H)$.

We denote by $P_m$ the $L_2$-orthogonal projection
onto the space $\mathrm{Span}\{w_1,\dots,w_m\}$, and we set $Q_m=I-P_m$.

\begin{definition}
We call a set of modes $\{w_j\}_{j=1}^m$
determining   (see  \cite{FMRT}, \cite{Foias-Prodi}) if
\begin{equation}\label{u-to-v}
\lim_{t\to\infty}\|u(t)-v(t)\|=0,
\end{equation}
as long as
$$
\lim_{t\to\infty}\|f(t)-g(t)\|=0
\quad\text{and}
\quad
\lim_{t\to\infty}\|P_m(u(t)-v(t))\|=0.
$$
Accordingly, a set of points $\{x^i\}_{i=1}^N\subset\Omega$ is called
a set of determining nodes (see~\cite{FMRT}, \cite{Foias-Temam})
if~(\ref{u-to-v}) holds as long as
$$
\lim_{t\to\infty}\|f(t)-g(t)\|=0
\quad\text{and}
\quad
\lim_{t\to\infty}\eta(u(t)-v(t))=0
$$
where $\eta(w)=\max_{j=1,\dots,N}|w(x^j)|$.
\end{definition}

We further suppose that
\begin{equation}\label{limsupfg}
\limsup_{t\to\infty}\|f(t)\|=:\mathbf{f}<\infty
\end{equation}

Subtracting in~(\ref{two-NS-eq}) the second equation from
the first and setting $w(t)=u(t)-v(t)$ and $h(t)=f(t)-g(t)$
we obtain
$$
\partial_tw+\nu Aw+B(u,w)+B(w,u)-B(w,w)=h(t).
$$
We write $w=p+q$, where $p=P_mw$, $q=Q_mw$ and
 take the scalar product  with $q$:
$$
\frac12\partial_t\|q\|^2+\nu\|\nabla\,q\|^2+
b(q,u,q)=(h,q)-b(u,p,q)-b(p,u,q)+b(p,p,q)+b(q,p,q),
$$
where $b(u,v,w)=(B(u,v),w)$.

A variant of the Gronwall lemma \cite{FMRT},
\cite{F-M-T-T}, \cite{Jones-Titi2}
is essential in the estimates below.
\begin{lemma}\label{L:Jones-Titi}
Suppose that $\alpha(t)$ and $\beta(t)$ are
 locally integrable  functions on $(0,\infty)$
satisfying for some $T>0$ the following conditions:
$$
\liminf_{t\to\infty}\frac1T\int\nolimits_t^{t+T}\alpha(\tau)d\tau=\gamma,
\
\limsup_{t\to\infty}\frac1T\int\nolimits_t^{t+T}\alpha^-(\tau)d\tau=
\Gamma,\
\lim_{t\to\infty}\frac1T\int\nolimits_t^{t+T}\beta^+(\tau)d\tau=0,
$$
where $\gamma>0$, $\Gamma<\infty$ and $\alpha^-=\max\{-\alpha,0\}$,
$\beta^+=\max\{\beta,0\}$. If $\xi(t)\ge0$, an
absolutely continuous function,  satisfies
$$
\xi'+\alpha\xi\le\beta\quad\text{on}\quad (0,\infty),
$$
then $\xi(t)\to0$ as $t\to\infty$.
\end{lemma}

All the terms on the right-hand side containing $h$ or $p$
can be absorbed in the function $\beta(t)$ in  Lemma~\ref{L:Jones-Titi}.
For example, using the second inequality in~(\ref{est-for-B})
we obtain
$$
\begin{aligned}
|b(u,p,q)|&\le c_2\|u\|^{1/2}\|\nabla\,u\|^{1/2}\|\nabla\,p\|
\|q\|^{1/2}\|\nabla\,q\|^{1/2}\\
&\le c_2(\lambda_m/\lambda_1)^{1/2}
\|p\|\|\nabla\,u\|\|\nabla\,q\|\le
(c_2/2)(\lambda_m/\lambda_1)^{1/2}\|p\|
(\|\nabla\,u\|^2+\|\nabla\,q\|^2).
\end{aligned}
$$
Our claim follows since $\|p(t)\|\to0$ and in view of~(\ref{time-average})
$$
\begin{aligned}
&\frac1T\int_{t}^{t+T}\|p(\tau)\|
(\|\nabla\,u(\tau)\|^2+\|\nabla\,q(\tau)\|^2)d\tau\\
\le
&\max_{\tau\in[t,t+T]}\|p(\tau)\|\
\frac1T\int_{t}^{t+T}
(\|\nabla\,u(\tau)\|^2+\|\nabla\,q(\tau)\|^2)d\tau\to0\quad\text{as}
\quad t\to\infty.
\end{aligned}
$$
The remaining terms can be treated in exactly the same way.
Therefore we obtain
$$
\partial_t\|q\|^2+2\nu\|\nabla\,q\|^2+
2b(q,u,q)\le\beta(t),
$$
where $\beta(t)\to0$ as $t\to\infty$.

Using the first inequality in~(\ref{est-for-B}) we  have
$$
2|b(q,u,q)|\le \nu\|\nabla\,q\|^2+c_b^2\nu^{-1}\|q\|^2\|\nabla\,u\|^2
$$
and then using the Poincar\'e inequality
$\lambda_{m+1}\|q\|^2\le\|\nabla\,q\|^2$ we obtain
$$
\partial_t\|q\|^2+\alpha(t)\|q\|^2\le\beta(t),\quad
\text{where}\quad
\alpha(t)=\nu\lambda_{m+1}-\nu^{-1}c_b^2\|\nabla\,u(t)\|^2.
$$
Since  by the well-known estimates for the Navier--Stokes system
(see, for instance,~\cite{CF88},\cite{T})
\begin{equation}\label{time-average}
\limsup_{t\to\infty}
\frac1T\int_t^{t+T}\|\nabla\,u(\tau)\|^2d\tau\le
\frac{\mathbf{f}^2}{T\nu^3\lambda_1^2}+
\frac{\mathbf{f}^2}{\nu^2\lambda_1}\,,
\end{equation}
it follows that
$\alpha$ satisfies the conditions of Lemma~\ref{L:Jones-Titi}
provided that $T$ is sufficiently large and
\begin{equation}\label{cond-lambda}
\lambda_{m+1}\,>\,
\frac{c_b^2 \mathbf{f}^2}{\nu^4\lambda_1}\,.
\end{equation}

In view of~(\ref{lower-bounds}) and~(\ref{est-for-B}), (\ref{est-for-cb})
 we have proved the following theorem.
\begin{theorem}\label{Thm-Dir}
The first $m$ eigenfunctions of the Stokes operator
are determining for the two-dimensional
Navier--Stokes system with Dirichlet boundary
conditions
if
\begin{equation}\label{cond-Dir}
m+1>\frac4{27\pi^3}\,G^2,\quad\text{where}\quad
G=\frac{\mathbf{f}|\Omega|}{\nu^2}\,.
\end{equation}
\end{theorem}\begin{remark}\label{Rem-bebore}
{\rm
The above theorem without explicit value of the constant
has been mentioned as a remark in~\cite{H-Ti}
and was also proved in~\cite{FMRT} and~\cite{Ju}.
}
\end{remark}

\subsection*{Periodic boundary conditions. Determining modes}
We now consider the Navier--Stokes system~(\ref{NSE})
with space periodic boundary conditions
$x\in \Omega=[0,L_1]\times[0,L_2]$.
We set $L_2=L$ and $L_1=L/\gamma$. Without
loss of generality we assume that $\gamma\le1$.
As before, $|\Omega|$ denotes the measure of the periodic domain
$\Omega$: $|\Omega|=L_1L_2=L^2/\gamma$.

 We further assume that
$u$, $p$ and $f$ have mean value zero over the torus.
Applying the  $\rot$ (curl) operator to the first equation in~(\ref{NSE})
we obtain the well-known vorticity equation
\begin{equation}\label{vort-eq}
\partial_t\vphi+J(\Delta^{-1}\vphi,\vphi)-\nu\Delta\,\vphi=\rot f,
\end{equation}
where  $\rot u=\vphi$, $u=\nabla^\perp\Delta^{-1}\vphi$,
$J(a,b)=\partial_1a\partial_2b-\partial_2a\partial_1b=
\nabla^\perp\,a\cdot\nabla\,b$,
$\nabla^\perp\psi=\mathbf{k}\times\nabla\,\psi=
(-\partial_2\psi,\partial_1\psi)$, and $\mathbf{k}$ is the
vertical unit vector.

We now recall 
 that the
spectrum~$\{\lambda_j\}_{j=1}^\infty$ of the
Stokes operator with periodic boundary conditions
coincides with that of the negative scalar Laplacian
$$
-\Delta\,\vphi_j=\lambda_j\vphi_j
$$
and the corresponding eigenfunctions are as follows:
$$
A\,w_j=\lambda_jw_j,\quad
w_j=\lambda_j^{-1/2}\nabla^\perp\vphi_j=
\lambda_j^{-1/2}(-\partial_2\vphi_j,\partial_1\vphi_j).
$$
Therefore the modes $\{w_1,\dots,w_m\}$ are determining for
the Navier--Stokes system~(\ref{NSE}) with periodic
boundary conditions  if the modes~$\{\vphi_1,\dots,\vphi_m\}$
are determining for the equation~(\ref{vort-eq}).

Similarly to~(\ref{two-NS-eq})
 we write
\begin{equation}\label{two-vort-eq}
\begin{aligned}
&\partial_t\vphi+J(\Delta^{-1}\vphi,\vphi)-\nu\Delta\,\vphi=\rot f(t),\\
&\partial_t\psi+J(\Delta^{-1}\psi,\psi)-\nu\Delta\,\psi=\rot g(t).
\end{aligned}
\end{equation}
Setting $\omega=\vphi-\psi$ and $H(t)=\rot f(t)-\rot g(t)$ we
obtain for $\omega$ the equation
\begin{equation}\label{eq-omega}
\partial_t\omega-\nu\Delta\,\omega+J(\Delta^{-1}\vphi,\omega)
+J(\Delta^{-1}\omega,\vphi)-J(\Delta^{-1}\omega,\omega)=H.
\end{equation}
As before, we write $\omega=p+q$, where $p=P_m\omega$ and
$q_m=Q_m\omega$ and where $P_m$ is the orthogonal projection
$P_m:L_2(\Omega)\to\mathrm{Span}(\vphi_1,\dots,\vphi_m)$.
Taking the scalar product with $q$ we obtain
$$
\aligned
&\frac12\partial_t\|q\|^2+\nu\|\nabla\,q\|^2+
(J(\Delta^{-1}q,\vphi),q)\\&=(H,q)-(J(\Delta^{-1}\vphi,p),q)-
(J(\Delta^{-1}p,\vphi),q)+(J(\Delta^{-1}p,p),q)+(J(\Delta^{-1}q,p),q)
=:\beta(t).
\endaligned
$$
As before, all the terms on the right-hand side containing $p$ can be
absorbed in $\beta(t)$ in Lemma~\ref{L:Jones-Titi}. Next we have
$$
\aligned
|(J(\Delta^{-1}q,\vphi),q)|&\le
\int|\nabla\Delta^{-1}q||\nabla \vphi|\,|q|\,dx
\le \|\nabla\Delta^{-1}q\|_{L_4}\|\nabla \vphi\|\|q\|_{L_4}\\
&\le c_\mathrm{L}(\gamma)^2\|\nabla\Delta^{-1}q\|^{1/2}
\|q\|\|\nabla\,q\|^{1/2}\|\nabla\vphi\|
\le\lambda_m^{-1/2}
 c_\mathrm{L}(\gamma)^2\|q\|\|\nabla\,\vphi\|\|\nabla\,q\|,
\endaligned
$$
where we used the Ladyzhenskaya inequality (see, for instance, \cite{CF88},
\cite{Lad}, \cite{TNS})
$$
\|\vphi\|_{L_4}\le c_\mathrm{L}(\gamma)
\|\vphi\|^{1/2}\|\nabla\,\vphi\|^{1/2},
\quad
\|\nabla\,\vphi\|_{L_4}\le c_\mathrm{L}(\gamma)\|\nabla\,\vphi\|^{1/2}
\|\Delta\,\vphi\|^{1/2}.
$$
We arrive at the same estimate if we use the
integral identity
$$
(J(f,g),h)=(J(h,f),g)
$$
and the Agmon inequality (see Theorem~\ref{Thm-cAT})
$$
\|\vphi\|_\infty\le c_{\mathrm{AT}}(\gamma)
\|\vphi\|^{1/2}\|\Delta\,\vphi\|^{1/2}.
$$
In fact,
$$
\aligned
|(J(\Delta^{-1}q,\vphi),q)|&=|(J(\vphi,q),\Delta^{-1}q)|
\le \|\Delta^{-1}q\|_\infty\|\nabla\, \vphi\|\|\nabla\,q\|\\
&\le c_\mathrm{AT}(\gamma)\|\Delta^{-1}q\|^{1/2}
\|q\|^{1/2}\|\nabla\,q\|\|\nabla\,\vphi\|
\le\lambda_m^{-1/2} c_\mathrm{AT}(\gamma)
\|q\|\|\nabla\,\vphi\|\|\nabla\,q\|,
\endaligned
$$
which gives
$$
|(J(\Delta^{-1}q,\vphi),q)|
\le\lambda_m^{-1/2} c_\mathrm{J}\|q\|\|\nabla\,\vphi\|\|\nabla\,q\|,
\quad\text{where}\quad
c_\mathrm{J}=\min(c_\mathrm{L}(\gamma)^2,c_\mathrm{AT}(\gamma)).
$$

As before we obtain the differential inequality
$$
\partial_t\|q\|^2+\alpha(t)\|q\|^2\le 2\beta(t),
\quad\text{where}\quad
\alpha(t)=\nu\lambda_{m+1}-
\nu^{-1}\lambda_{m+1}^{-1}c_\mathrm{J}^2\|\nabla\,\vphi(t)\|^2.
$$
It follows from the well-known a priory estimate on the time average of the $H^2$-norm
of a solution $u$ (see, for instance, \cite{BV}, \cite{CF88},
\cite{Jones-Titi2}, \cite{T})
\begin{equation}\label{apriori-est}
\limsup_{t\to\infty}
\frac1T\int_t^{t+T}\|\nabla\,\vphi(\tau)\|^2d\tau=
\limsup_{t\to\infty}
\frac1T\int_t^{t+T}\|A\,u(\tau)\|^2d\tau\le
\frac{\mathbf{f}^2}{T\nu^3\lambda_1}+
\frac{\mathbf{f}^2}{\nu^2}\,
\end{equation}
 that
$\alpha$ satisfies conditions of Lemma~\ref{L:Jones-Titi}
provided that $T$ is sufficiently large and
\begin{equation}\label{cond-lambda-per}
\lambda_{m+1}^2\,>\,
\frac{c_\mathrm{J}^2 \mathbf{f}^2}{\nu^4}\,.
\end{equation}
It was shown in~\cite{I-M-T}, \cite{I05MS}
 that $c_\mathrm{L}(\gamma)\le(6/(\gamma\pi))^{1/4}$.
In the Appendix in section~\ref{S:Appen} we will show that
$c_\mathrm{AT}(\gamma)\le1/\sqrt{\gamma\pi}$.
Hence we can take $c_\mathrm{J}=1/\sqrt{\gamma\pi}$.
Furthermore, for $\gamma=1$
$\lambda_m\ge(\lambda_1/4)m$, where $\lambda_1=4\pi^2L^{-2}$.
We obtain the following theorem.
\begin{theorem}\label{Thm-per}
The first $m$ eigenfunctions of the Stokes
operator are determining for the two-dimensional
Navier--Stokes system with periodic boundary
conditions if
\begin{equation}\label{cond-per-gamma}
\lambda_{m+1}>
\left(\frac1{\gamma\pi}\right)^{1/2}\frac{\mathbf{f}}{\nu^2}.
\end{equation}

For a square torus $(\gamma=1)$ this condition
is satisfied if
\begin{equation}\label{cond-per}
m+1>
\frac1{\pi^{3/2}}\,G,\quad\text{where}\quad
G=\frac{\mathbf{f}L^2}{\nu^2}\,.
\end{equation}
\end{theorem}

\begin{remark}\label{Rm-lb-gamma}
{\rm
The first eigenvalues $\lambda_1$, $\lambda_2$, $\dots$
of the Laplacian on the periodic domain
$\Omega=[0,L/\gamma]\times[0,L]$
are of order $\gamma^2$ when $\gamma\ll 1$.
It was shown in~\cite{I05MS} (see Proposition~4.1) that if
$m\ge2/\gamma$, then
\begin{equation}\label{lb-gamma}
\lambda_m\ge\frac{m\gamma}8\cdot\frac{4\pi^2}{L^2}=
\frac{\pi^2}2\cdot\frac m{|\Omega|}\,.
\end{equation}
Therefore condition~(\ref{cond-per-gamma}) is satisfied if
\begin{equation}\label{cond-per-gamma-expl}
m+1>\frac2\gamma\,+\,
\frac2{\pi^2}\left(\frac1{\gamma\pi}\right)^{1/2}
\frac{\mathbf{f}|\Omega|}{\nu^2}\,.
\end{equation}

}
\end{remark}

\subsection*{Determining nodes}
Suppose that the periodic domain $\Omega$ is divided into
$N$ equal squares with side of length $l$ and for each square there
is a point $x^j$, $j=1,\dots,N$ chosen arbitrarily in it.
For the two solutions $u$ and $v$ of the Navier-Stokes
equations we assume that
$$
\eta(w(t))=\max_{j=1,\dots,N}|w(t,x^j)|\to0\quad\text{as}
\quad t\to\infty,
\quad\text{where}\quad
w=u-v.
$$

We take the scalar product of~(\ref{eq-omega}) and $\omega$:
$$
\partial_t\|\omega\|^2+2\nu\|\nabla\,\omega\|^2=
\beta(t)-2(J(\Delta^{-1}\omega,\vphi),\omega),
\quad\text{where}\quad\beta(t)=2(H,\omega).
$$
For the nonlinear term using inequalities~(\ref{Ag-Tai-vec}),
(\ref{in-vort-form}) and Young's inequality
 we have
$$
\aligned
2|(J(\Delta^{-1}\omega,\vphi),\omega)|\le
2\|\nabla\Delta^{-1}\omega\|_\infty\|\nabla\,\vphi\|\|\omega\|\le
2c_\mathrm{AT}(\gamma)\|\nabla\Delta^{-1}\omega\|^{1/2}\|\nabla\,\omega\|^{1/2}
\|\omega\|\|\nabla\,\vphi\|\le\\
2c_\mathrm{AT}(\gamma)
(4|\Omega|)^{1/4}\eta^{1/2}(w)
\|\nabla\,\omega\|^{1/2}\|\omega\|\|\nabla\,\vphi\|+
2c_\mathrm{AT}(\gamma)68^{1/4}l\|\nabla\,\omega\|
\|\omega\|\|\nabla\,\vphi\|\le\\
\beta_1(t)+
c_\mathrm{AT}(\gamma)^2
68^{1/2}l^2\nu^{-1}\|\omega\|^2\|\nabla\,\vphi\|^2+
\nu\|\nabla\,\omega\|^2.
\endaligned
$$
Hence
$$
\partial_t\|\omega\|^2+\nu\|\nabla\,\omega\|^2\le
\beta_1(t)+
c_\mathrm{AT}(\gamma)^2 68^{1/2}l^2\nu^{-1}\|\omega\|^2\|\nabla\,\vphi\|^2.
$$
Using inequality~(\ref{Jones-Titi2})
in the form
$\|\nabla\,\omega\|^2\ge 68^{-1/2}l^{-2}\|\omega\|^2-\beta_2(t)$
to bound from below the second term
on the left-hand side we obtain
$$
\partial_t\|\omega\|^2+\alpha(t)\|\omega\|^2\le
\beta(t),
$$
where
$$
\alpha(t)=68^{-1/2}\nu l^{-2}-c_\mathrm{AT}(\gamma)^2 68^{1/2}\nu^{-1}l^2
\|\nabla\,\vphi\|^2.
$$
 Taking into account~(\ref{apriori-est})
we see that if the number of squares
 $N=|\Omega|/l^2$ is sufficiently large
 (or, equivalently, the typical distance $l$ between the nodes
 is sufficiently small),
then $\alpha$ satisfies conditions of Lemma~\ref{L:Jones-Titi} and
 the corresponding $N$ nodes are determining.
We obtained the following result.
\begin{theorem}\label{T:nodes-NS} If
\begin{equation}\label{det-nodes-NS}
N>\left(\frac{68}{\gamma\pi}\right)^{1/2}\frac{\mathbf{f}|\Omega|}{\nu^2}
=
\left(\frac{68}{\gamma^3\pi}\right)^{1/2}\frac{\mathbf{f}L^2}{\nu^2}
\end{equation}
equal squares tile $\Omega=[0,L/\gamma]\times[0,L]$,
then $N$ nodes $($chosen arbitrarily one in each square$)$
are determining for the
space periodic Navier--Stokes system in $\Omega$.
\end{theorem}

\begin{remark}\label{Rem-bebore-per}
{\rm
 The estimates for determining modes and nodes of
 Theorems~\ref{Thm-per} and~\ref{T:nodes-NS}
were obtained  for $\gamma=1$ in~\cite{Jones-Titi2}.
}
\end{remark}

\setcounter{equation}{0}
\section{Determining modes and nodes for damped
Navier--Stokes equations}\label{S:damped N-S}

In this section we consider the damped-driven Navier--Stokes
system having important applications in geophysical
hydrodynamics~\cite{D-F}, \cite{P}
\begin{equation}\label{NSE-damped}
\begin{aligned}
\partial_tu+\sum_{i=1}^2u^i\partial_iu&=
\nu\Delta\,u-\mu u-\nabla\,p +f,\\
\div u&=0,\quad
u(0)=u_0.
\end{aligned}
\end{equation}
\subsection*{Periodic boundary conditions}
We first consider this system on the torus
 $x\in \Omega=[0,L/\gamma]\times[0,L]$, with space periodic
boundary conditions. The right-hand side $f=f(t)$
satisfies the condition
\begin{equation}\label{cond-for-f}
    \limsup_{t\to\infty}\|\rot f(t)\|_\infty=:\mathbf{F}_\infty<\infty.
\end{equation}
\begin{lemma}\label{L:Est-for-rot}
The following bound holds for $u(t)$:
\begin{equation}\label{rot-infty}
    \limsup_{t\to\infty}\|\rot u(t)\|_\infty\le
    \frac{\mathbf{F}_\infty}{\mu}\,.
\end{equation}
\end{lemma}
\begin{proof}
We use the vorticity formulation of the system~(\ref{NSE-damped})
\begin{equation}\label{vort-form}
\partial_t\vphi+J(\Delta^{-1}\vphi,\vphi)-\nu\Delta\,\vphi+\mu\vphi=\rot f
\end{equation}
and take the scalar product with $\vphi^{2k-1}$, where $k\ge1$ is
integer,
and use the identity $(J(\psi,\vphi),\vphi^{2k-1})=
(2k)^{-1}\int J(\psi,\vphi^{2k})dx=
(2k)^{-1}\int\div(\vphi^{2k}\nabla^\perp\psi)dx=0$.
We obtain
$$
\aligned
\|\vphi\|_{L_{2k}}^{2k-1}\partial_t\|\vphi\|_{L_{2k}}+
(2k-1)\nu\int|\nabla\,\vphi|^2\vphi^{2k-2}dx+
\mu\|\vphi\|_{L_{2k}}^{2k}=\\=(\rot f(t),\vphi^{2k-1})\le\
\|\rot f(t)\|_{L_{2k}}\|\vphi\|_{L_{2k}}^{2k-1}.
\endaligned
$$
Hence, by Gronwall's inequality
$$
\|\vphi(t+\tau)\|_{L_{2k}}\le\|\vphi(\tau)\|_{L_{2k}}e^{-\mu t}+
\mu^{-1}\sup_{s\in[\tau,\infty)}\|\rot f(s)\|_{L_{2k}}(1-e^{-\mu t}),
$$
and passing to the limit as $k\to\infty$ we find
$$
\|\vphi(t+\tau)\|_\infty\le\|\vphi(\tau)\|_\infty e^{-\mu t}+
\mu^{-1}\sup_{s\in[\tau,\infty)}\|\rot f(s)\|_\infty(1-e^{-\mu t}).
$$
Now, we let $t\to\infty$ to obtain
$$
\limsup_{t\to\infty}\|\vphi(t)\|_\infty\le\frac{\mathbf{F}_\infty}\mu\,.
$$
\end{proof}

We consider the systems~(\ref{NSE-damped})
with right-hand sides $f$ and $g$ such that
$$
\lim_{t\to\infty}\|\rot (f(t)-g(t))\|_\infty=0.
$$
Similarly to (\ref{two-vort-eq}) for
$\omega=\vphi-\psi$  we
obtain the equation
\begin{equation}\label{eq-for-q-damped}
\partial_t\omega-\nu\Delta\,\omega+\mu \omega+J(\Delta^{-1}\vphi,\omega)
+J(\Delta^{-1}\omega,\vphi)-J(\Delta^{-1}\omega,\omega)=H,
\end{equation}
 where $H(t)=\rot f(t)-\rot g(t)$.

\subsection*{Determining modes}
We take the scalar product of~(\ref{eq-for-q-damped}) with
$q=Q_m\omega$:
$$
\aligned
\partial_t\|q\|^2+2\nu\|\nabla\,q\|^2+2\mu\|q\|^2
&\le 2\beta(t)+2|(J(\Delta^{-1}q,\vphi),q)|\le\\
2\beta(t)+2\|\nabla q\|\|\nabla\Delta^{-1}q\|\|\vphi\|_\infty&\le
2\beta(t)+
\nu\|\nabla\,q\|^2+\nu^{-1}\lambda_{m+1}^{-1}\|q\|^2\|\vphi\|_\infty^2,
\endaligned
$$
where we used Young's and Poincar\'e inequalities.
Dropping the $\mu$-term on the left-hand side and again
using the Poincar\'e inequality we obtain
$$
\partial_t\|q\|^2+\|q\|^2
(\nu\lambda_{m+1}-\nu^{-1}\lambda_{m+1}^{-1}\|\vphi(t)\|_\infty^2)\le
2\beta(t).
$$
By estimate~(\ref{rot-infty}) and  Lemma~\ref{L:Jones-Titi} the
first $m$ modes are determining provided
\begin{equation}\label{cond-modes-damped}
\lambda_{m+1}\ge\frac{\mathbf{F}_\infty}{\mu\nu}\,.
\end{equation}
Using~(\ref{lb-gamma}) we see that this condition is satisfied if
\begin{equation}\label{cond-modes-damped-expl}
m+1>\max\left\{\frac2\gamma\,,\
\frac 2{\pi^2}\frac{\mathbf{F}_\infty|\Omega|}{\mu\nu}
\right\}.
\end{equation}

For a square torus $\lambda_m\ge\lambda_1/4$, $\lambda_1=4\pi^2/L^2$,
hence condition~(\ref{cond-modes-damped}) is satisfied if
\begin{equation}\label{m-cond-modes-damped}
m+1\ge\frac1{\pi^2}\frac{\mathbf{F}_\infty L^2}{\mu\nu}\,.
\end{equation}

\subsection*{Determining nodes}
Suppose that our periodic domain $\Omega$ is divided into $N$
equal squares $Q_j$ with side of length $l$, $j=1,\dots,N$ and
we chose arbitrarily a point $x^j\in Q_j$
for each $j=1,\dots,N$.

For $u\in H^2_\mathrm{per}(\Omega)$ we set
\begin{equation}\label{eta}
    \eta(u)=\max_{j=1,\dots,N}|u(x^j)|.
\end{equation}
Suppose that $\eta(w(t))\to0$ as $t\to\infty$ for
$w(t)=u(t)-v(t)$, where $u$ and $w$ are two solutions
of~(\ref{NSE-damped}).

We take the scalar product
of~(\ref{eq-for-q-damped}) with $\omega$:
$$
\frac12\partial_t\|\omega\|^2+\nu\|\nabla\,\omega\|^2+\mu\|\omega\|^2=
(H,\omega)-(J(\Delta^{-1}\omega,\vphi),\omega).
$$
We estimate
the nonlinear term by means of inequality~(\ref{in-vort-form})
$$
\aligned
|(J(\Delta^{-1}\omega,\vphi),\omega)|=
|(J(\omega,\Delta^{-1}\omega),\vphi)|\le
\|\nabla\,\omega\|\|\nabla\Delta^{-1}\omega\|\|\vphi\|_\infty\le\\
2|\Omega|^{1/2}\eta(w)\|\nabla\,\omega\|\|\vphi\|_\infty+
\sqrt{68}|\Omega|N^{-1}\|\nabla\,\omega\|^2\|\vphi\|_\infty\,.
\endaligned
$$
As a result we obtain
$$
\frac12\partial_t\|\omega\|^2+\alpha(t)\|\nabla\,\omega\|^2
\le\beta(t),
$$
where
$$
\alpha(t)=\nu-\sqrt{68}|\Omega|N^{-1}\|\vphi(t)\|_\infty,
\quad\
\beta(t)=(H(t),\omega)+
2|\Omega|^{1/2}\eta(w(t))\|\nabla\,\omega\|\|\vphi(t)\|_\infty\,.
$$
As before $\beta(t)\to0$ as $t\to\infty$,
while in view of~(\ref{rot-infty}) $\alpha(t)\ge\mathrm{const}>0$
for all $t$ large enough provided that
$N>\sqrt{68}\,{ \mathbf{F}_\infty|\Omega|}/({\mu\nu})\,.$

We combine  the results so obtained above in the following theorem.

\begin{theorem}\label{Th:modes-nodes}
The first $m$ modes of the Stokes operator are determining
for the Navier--Stokes system with damping~$(\ref{NSE-damped})$ if
\begin{equation}\label{cond-modes-damped-expl+}
m+1>\frac2\gamma\,+\,
\frac 2{\pi^2}\frac{\mathbf{F}_\infty|\Omega|}{\mu\nu}.
\end{equation}

If $\Omega$ is tiled by $N$ equal squares, then any collection
of nodes $($one in each square$)$ is determining if
\begin{equation}\label{cond-for-nodes}
N>\sqrt{68}\,\frac{ \mathbf{F}_\infty|\Omega|}{\mu\nu}\,.
\end{equation}
\end{theorem}
\begin{remark}\label{Rem:aspect-ratio}
{\rm
We observe that estimates for the number of the determining modes and
nodes~(\ref{cond-modes-damped-expl}) and (\ref{cond-for-nodes})
depend linearly on the measure of the
periodic domain $|\Omega|$ and depend on
 the aspect ratio $\gamma$ of the torus
only via $|\Omega|=L^2/\gamma$. Furthermore,  the characteristic
microscopic length $l$ of the lattice of the determining nodes
 satisfies the following $\Omega$-independent
estimate from above
\begin{equation}\label{est-for-l}
l<68^{-1/4}\left(
\frac{\mu\nu}{\mathbf{F}_\infty}
\right)^{1/2}.
\end{equation}
}
\end{remark}

\begin{remark}\label{Rem:fractal}
{\rm
It was shown in~\cite{I-M-T} that the fractal dimension of the global
attractor of the autonomous system~(\ref{NSE-damped}) on the
torus $[0,L]^2$ satisfies the estimate
\begin{equation}\label{dim-upper}
\dim_F\mathcal{A}\le\left(\frac 6{\pi^3}\right)^{1/2}
\frac{\|\rot f\|L}{\mu\nu}\,.
\end{equation}
It was also shown that for the Kolmogorov forcing
of the form
\begin{equation}\label{Kolmf}
f=f_s=\begin{cases}f_1=
A(\mu,\nu)\sin s \frac{2\pi x_2}L,\\
f_2=0,\end{cases}
\end{equation}
one has a lower bound
\begin{equation}\label{dim-lower}
\dim_F\mathcal{A}\ge\mathrm{const}
\frac{\|\rot f\|L}{\mu\nu}
\end{equation}
in the limit $\nu\to0^+$ (accordingly, $s=(\mu L^2/\nu)^{1/2}\to\infty$).

Since $\|\rot f\|\le\|\rot f\|_\infty L$,
it follows that
\begin{equation}\label{dim-upper-infty}
\dim_F\mathcal{A}\le\left(\frac 6{\pi^3}\right)^{1/2}
\frac{\|\rot f\|_\infty L^2}{\mu\nu}\,,
\end{equation}
 and since
$\|\rot f_s\|=(L/\sqrt{2})\|\rot f_s\|_\infty$,
it follows from~(\ref{dim-lower}) that the estimate~(\ref{dim-upper-infty})
 is also sharp
with respect to the dimensionless number
$\frac{\|\rot f\|_\infty L^2}{\mu\nu}$.

Hence the  bounds~(\ref{m-cond-modes-damped}) and~(\ref{cond-for-nodes})
for the number of determining modes and nodes
for the damped Navier--Stokes system
are of the  same order as the fractal dimension of the global attractor.

We point out that in general and for the Navier--Stokes equations
there is a gap between the number of the determining modes and nodes and
the dimension of the global attractor. However, the works
\cite{F-O}, \cite{H-K} and \cite{Friz-Robinson}
indicate that one can perturb the points   or the projections
to obtain the number of nodes and the rank of the projections
comparable with the dimension of the global attractor.
Here, however, we show that there is no need for perturbation
and that the usual projections $P_m$  and any choice of points
will do. It will therefore be interesting to
understand the role of damping term here in terms
of the generic results of \cite{F-O}, \cite{H-K}, \cite{Friz-Robinson}
which rely heavily on the  M\~an\'e embedding theorem.
}
\end{remark}

\subsection*{Stress-free boundary conditions}
Let $\Omega\subset\mathbb{R}^2$ be a bounded simply connected
domain with $C^2$ boundary. Let $n$
be the outward unit normal vector.
 We consider the system~(\ref{NSE-damped})
supplemented with the so-called stress-free boundary conditions
\begin{equation}\label{stress-free}
u\cdot n|_{\partial\Omega}=0,
\qquad
\rot u|_{\partial\Omega}=0.
\end{equation}
Then any smooth vector field $u$, $\div u=0$,
satisfying~(\ref{stress-free}) has a unique single valued
stream function $\psi$, $u=\nabla^\perp\psi$
with $\psi|_{\partial\Omega}=0$ and
$\Delta\,\psi|_{\partial\Omega}=\rot u|_{\partial\Omega}=0$.
Therefore, the vorticity formulation for the system~(\ref{NSE-damped}),
(\ref{stress-free}) is the equation~(\ref{vort-form})
with zero boundary condition both for $\vphi$ and $\psi$:
\begin{equation}\label{vort-form-mu}
\aligned
\partial_t\vphi+J(\psi,\vphi)&-
\nu\Delta\,\vphi+\mu\vphi=\rot f,\\
\Delta\,\psi&=\vphi,\\
\vphi|_{\partial\Omega}&=
\psi|_{\partial\Omega}=0.
\endaligned
\end{equation}

For the Stokes eigenvalue problem with boundary
conditions~(\ref{stress-free})
\begin{equation}\label{Stokes-sterss-free}
\aligned
-\Delta\, w_k+\nabla\,p_k=\lambda_kw_k,
\quad \div w_k=0,\\
w_k\cdot n|_{\partial\Omega}=0, \quad \rot w_k|_{\partial\Omega}=0,
\endaligned
\end{equation}
we have (as in the case of periodic boundary conditions)
that
 $\{\lambda_k\}_{k=1}^\infty$
 are the eigenvalues of the scalar Dirichlet
problem $-\Delta\,\vphi_k=\lambda_k\vphi_k$,
$\vphi_k|_{\partial\Omega}=0$
and $w_k=\lambda_k^{-1/2}\nabla^\perp\vphi_k$.
Hence  we can use the Li--Yau lower bound \cite{Li-Yau}
 for the eigenvalues $\lambda_k$
\begin{equation}\label{Li-Yau}
\lambda_k\ge\frac{2\pi k}{|\Omega|}\,.
\end{equation}
Lemma~\ref{L:Est-for-rot} and the subsequent argument for estimates
of the determining modes still hold and we
obtain that condition~(\ref{cond-modes-damped}) is sufficient
for the first $m$ modes to be determining.
In view of~(\ref{Li-Yau}) we obtain the following result.

\begin{theorem}\label{Thm-stress-free}
The first $m$ modes of the Stokes operator are determining
for the Navier--Stokes system with damping~$(\ref{NSE-damped})$
with stress-free boundary conditions~$(\ref{stress-free})$ if
$$
m+1>\frac1{2\pi}\frac{ \mathbf{F}_\infty|\Omega|}{\mu\nu}\,.
$$
\end{theorem}

\begin{remark}\label{Rm-stress-free}
{\rm
As in the space-periodic case this estimate agrees with
the estimate for the fractal dimension of the attractor~\cite{I-M-T},
but, unlike the latter, does not involve constants depending
on the smoothness and shape of the boundary.
}
\end{remark}
\begin{remark}\label{Rm-stress-free-nodes}
{\rm
A similar result holds for the determining nodes
and other determining functionals and projections
 (see~\cite{Cockburn-Jones-Titi-1}, \cite{Cockburn-Jones-Titi})
if we use
extension operators mapping Sobolev spaces defined on $\Omega$
to the spaces defined on  corresponding periodic rectangular domain
containing $\Omega$. In this case, however,
the estimate involves a
constant depending
on the smoothness and shape of the boundary.
}
\end{remark}

\setcounter{equation}{0}
\section{Appendix. Proof of auxiliary inequalities}\label{S:Appen}

The embedding of the Sobolev space $H^l(M)$
with norm $\|u\|_{H^l}^2=\|u\|^2+\|(-\Delta)^{l/2}u\|^2$
into the space
of bounded continuous functions $C(M)$, where $\dim M=n$
and $l>n/2$, can be written as a multiplicative
inequality
\begin{equation}\label{Sob-to-C}
\|u\|_\infty\le c_M(l)\|u\|^{\theta}\|(-\Delta)^{l/2}u\|^{1-\theta},
\quad\text{where}\quad\theta=(2l-n)/2l.
\end{equation}
Inequalities of this type are sometimes called the Agmon
inequalities (see~\cite{Agmon}).
The best constant $c_M(l)$ for $M=\mathbb{R}$ in this inequality
was found in~\cite{Taikov} (the results of~\cite{Taikov} can easily
be generalized to the case when $M=\mathbb{R}^n$). Sharp constants
in inequalities for periodic functions and functions defined on
the sphere were found in~\cite{I98MS}. Following~\cite{I98MS} we
 consider below the case of a
two-dimensional torus.

The constant in inequality~(\ref{Sob-to-C}) on a two-dimensional
torus clearly depends only on the aspect ratio
$\gamma$ of the torus. We first consider the case of a
square torus $\gamma=1$, and then
without loss of generality we assume that
$\Omega=\mathrm{T}^2=[0,2\pi]^2$.

We consider the negative Laplacian  $-\Delta$ in
$H=L_2(\mathrm{T}^2)\cap\{\vphi,\ \int\vphi\,dx=0\}$ and order
its eigenvalues according to magnitude and
multiplicity:
\begin{equation}\label{alleig}
    1=\lambda_1\le\lambda_2\le\dots,
\qquad\{\lambda_j,\ j=1,\dots\}=\{k^2=k_1^2+k_2^2,\ k=(k_1,k_2)\in
\mathbb{Z}^2_0\},
\end{equation}
where $\mathbb{Z}^2_0=\mathbb{Z}^2\setminus\{0\}$.
 The corresponding basis of
orthonormal eigenfunctions $w_j(x)$, $-\Delta\,w_j=\lambda_jw_j$, is
the basis of trigonometric functions
\begin{equation}\label{basis}
\aligned
\bigcup_{j\in\mathbb{N}}w_j(x)= \bigcup_{k\in\mathbb{Z}^2_+}
\left\{({\sqrt{2}\pi})^{-1}\sin kx,\quad({\sqrt{2}\pi})^{-1}\cos kx\right\},\\
\Zbb_+^2=\{k\in\Zbb^2_0,\quad k_1\ge0,\ k_2\ge 0\}
\cup\{k\in\Zbb^2_0,\quad k_1\ge1,\ k_2\le 0\}.
\endaligned
\end{equation}
Similarly to~(\ref{alleig}), we write
\begin{equation}\label{halfeig}
    1=\Lambda_1\le\Lambda_2\le\dots,
\qquad\{\Lambda_j,\ j=1,\dots\}=\{k^2,\ k\in \mathbb{Z}^2_+\}
\end{equation}
and observe that
\begin{equation}\label{2lambda}
\bigcup_{j=1}^\infty\{\lambda_j\}=
\bigcup_{l=1}^\infty\{\Lambda_l,\Lambda_l\}.
\end{equation}
Hence, for $j\ge1$, we have $\Lambda_j=\lambda_{2j}=\lambda_{2j-1}$ and,
corresponding to each $\Lambda=\Lambda_j$, there are two eigenfunctions
$u_j(x)=(\sqrt{2}\pi)^{-1}\sin kx$ and $v_j(x)=(\sqrt{2}\pi)^{-1}\cos kx$
for some uniquely defined $k_j=(k_1(j), k_2(j))$ with $k_j^2=\Lambda_j$.
We obviously have
\begin{equation}\label{trig}
u_j(x)^2+v_j(x)^2=\frac1{2\pi^2}\,.
\end{equation}

\begin{theorem}\label{Thm-cAT}
The sharp constant $c_\mathrm{AT}$ in the inequality
\begin{equation}\label{Agmon-Taikov}
\|\vphi\|_\infty\,\le\,c_\mathrm{AT}\|\vphi\|^{1/2}\|\Delta\,\vphi\|^{1/2},
\quad\vphi\in H\cap H^2_\mathrm{per}(\mathrm{T}^2),
\end{equation}
is given by
\begin{equation}\label{const-cAT}
c_\mathrm{AT}^2=\frac1{\pi^2}\sup_{\mu>0}\,\mu\,\sum_{n=1}^\infty
\frac1{\mu^2+\Lambda_n^2}\,
\end{equation}
and, in particular,
\begin{equation}\label{const-cAT-expl}
c_\mathrm{AT}^2\,<\,\frac1{\pi}\,.
\end{equation}
\end{theorem}
\begin{proof}
Writing $\vphi$ in terms of the Fourier series
$\vphi(x)=\sum_{n=1}^\infty c_n w_n(x)$,
for an arbitrary point $x_0$ and a positive parameter $\nu$ we have
\begin{equation}\label{additive}
\begin{aligned}
\vphi(x_0)^2=\biggl(\sum_{n=1}^\infty c_n w_n(x_0)\biggr)^2\le
\sum_{n=1}^\infty
\frac{w_n(x_0)^2}{1+\nu\lambda_n^2}\
\sum_{n=1}^\infty c_n^2(1+\nu\lambda_n^2)=\\
\frac1{2\pi^2}
\sum_{n=1}^\infty
\frac{1}{1+\nu\Lambda_n^2}\cdot
\bigl(\|\vphi\|^2+\nu\|\Delta\,\vphi\|^2\bigr),
\end{aligned}
\end{equation}
where we used~(\ref{2lambda}), (\ref{trig}).
Since the right-hand side of~(\ref{additive}) is independent of~$x_0$,
it follows that
\begin{equation}\label{infty-additive}
\|\vphi\|_\infty^2\,\le\,
\frac1{2\pi^2}
\sum_{n=1}^\infty
\frac{1}{1+\nu\Lambda_n^2}\cdot
\bigl(\|\vphi\|^2+\nu\|\Delta\,\vphi\|^2\bigr).
\end{equation}
Let $x_0$ be fixed. Then there is equality in~(\ref{additive}),
(\ref{infty-additive}) if and only if
$$
c_n=(1+\nu\lambda_n^2)w_n(x_0),
$$
that is, if
\begin{equation}\label{cond-equal}
\vphi(x)=\sum_{n=1}^\infty
\frac{w_n(x)w_n(x_0)}{1+\nu\lambda_n^2}=
\frac{1}{2\pi^2}
\sum_{n=1}^\infty
\frac{\cos(k_n(x-x_0))}{1+\nu\Lambda_n^2}\,.
\end{equation}
We now set $\nu=\nu_*=\|\vphi\|^2/\|\Delta\,\vphi\|^2$. Then
$\|\vphi\|^2+\nu_*\|\Delta\,\vphi\|^2=
2\nu_*^{1/2}\|\vphi\|\|\Delta\,\vphi\|$ and therefore
$$
\|\vphi\|_\infty^2\,\le\,
\frac1{\pi^2}\,\nu_*^{1/2}
\sum_{n=1}^\infty
\frac{1}{1+\nu_*\Lambda_n^2}\cdot
\|\vphi\|\|\Delta\,\vphi\|\le
\frac1{\pi^2}\,\sup_{\nu>0}\,\nu^{1/2}\,
\sum_{n=1}^\infty
\frac{1}{1+\nu\Lambda_n^2}\cdot
\|\vphi\|\|\Delta\,\vphi\|,
$$
which shows (with $\nu=\mu^{-2}$) that $c_\mathrm{AT}^2$ is less than
or equal to the
right-hand side of~(\ref{const-cAT}).

Suppose now that the supremum of the function
$H(\nu)=\nu^{1/2}\sum_{n=1}^\infty(1+\nu\Lambda_n^2)^{-1}$
is attained at a finite point $\nu_*$, $0<\nu_*<\infty$. Then
\begin{equation}\label{two-equations}
\frac12\sum_{n=1}^\infty\frac1{1+\nu_*\Lambda_n^2}=
\nu_*\sum_{n=1}^\infty\frac{\Lambda_n^2}{(1+\nu_*\Lambda_n^2)^2}\,,
\qquad
\frac12\sum_{n=1}^\infty\frac1{1+\nu_*\Lambda_n^2}=
\sum_{n=1}^\infty\frac{1}{(1+\nu_*\Lambda_n^2)^2}\,.
\end{equation}
In fact, the first equality follows from $H'(\nu_*)=0$.
Summing the the first and the second equalities we
obtain a valid identity, hence the second equality
also holds.

Next we set $\nu=\nu_*$ in~(\ref{cond-equal}) and $x_0=0$. Then
for the corresponding $\vphi=\vphi_*$ we have
$$
\|\vphi_*\|^2=\sum_{n=1}^\infty\frac1{(1+\nu_*\Lambda_n^2)^2},
\qquad
\|\Delta\,\vphi_*\|^2=
\sum_{n=1}^\infty\frac{\Lambda_n^2}{(1+\nu_*\Lambda_n^2)^2}\,.
$$
Then it follows from~(\ref{two-equations}) that
 $\|\vphi_*\|^2/\|\Delta\,\vphi_*\|^2=\nu_*$.
Hence
$$
\|\vphi_*\|_\infty^2=\frac1{\pi^2}H(\nu_*)\|\vphi_*\|\|\Delta\,\vphi_*\|,
$$
which proves the theorem in the case when $0<\nu_*<\infty$.

Suppose now that supremum is attained as $\nu\to0$ (observe that
$H(\nu)\to0$ as $\nu\to\infty$). We consider
inequality~(\ref{Agmon-Taikov}) on the finite dimensional space
$T_N=\mathrm{Span}\{\sin kx,\ \cos kx\}$, $k^2\le N$.
The corresponding sharp constant
$c_\mathrm{AT}(N)$ is given by the formula
$$
c_\mathrm{AT}(N)^2=\frac1{\pi^2}\max_{\nu>0}\,H_N(\nu),
\quad H_N(\nu)=\nu^{1/2}
\sum_{\Lambda_n\le N}
\frac1{1+\nu\lambda_n^2}\,.
$$
The maximum is attained since $H_N(0)=0$ and $H_N(\nu)\to0$
as $\nu\to\infty$. Hence there exists an extremal function
$\vphi_*^N\in T_N$. Since the spaces $T_N$ are dense in the
Sobolev space $H^2_\mathrm{per}(\mathrm{T}^2)$, it follows that
$$
c_\mathrm{AT}=\lim_{N\to\infty}c_\mathrm{AT}(N),
\quad\text{and}\quad
c_\mathrm{AT}^2=\frac1{\pi^2}\sup_{\nu>0}\nu^{1/2}
\,\sum_{n=1}^\infty\frac1{1+\nu\Lambda_n^2}\,.
$$

It remains to prove~(\ref{const-cAT-expl}). Using the lower bound
$\Lambda_n\ge n/2$  (see~\cite{I-M-T}) we have
$$
\nu^{1/2} \,\sum_{n=1}^\infty\frac1{1+\nu\Lambda_n^2}=
\frac1\mu\sum_{n=1}^\infty\frac1{1+(\Lambda_n/\mu)^2}<
\frac1\mu\sum_{n=1}^\infty\frac1{1+(n/(2\mu))^2}<
\int_0^\infty f(x)\,dx=\pi,
$$
where $f(x)=1/(1+(x/2)^2)$ is monotone decreasing and the third
term in the above formula is the Riemann sum with step $1/\mu$
for the corresponding integral.
\end{proof}
\begin{remark}\label{Rem-vector}
{\rm
In the vector case
$u\in H^2_\mathrm{per}(\mathrm{T}^2)^2$ we have the same constant in the
corresponding inequality
\begin{equation}\label{Ag-Tai-vec}
\|u\|_\infty\le c_\mathrm{AT}\|u\|^{1/2}\|\Delta\, u\|^{1/2}.
\end{equation}
In fact, for $u=\{u^1,u^2\}$ we have
$$
\aligned
\|u\|_\infty^2\le\|u^1\|_\infty^2+\|u^2\|_\infty^2\le
c_\mathrm{AT}^2(\|u^1\|\|\Delta\,u^1\|+\|u^2\|\|\Delta\,u^2\|)\le\\
\frac{c_\mathrm{AT}^2}2(\ve(\|u^1\|^2+\|u^2\|^2)+
\ve^{-1}(\|\Delta\,u^1\|^2+\|\Delta\,u^2\|^2))=
\frac{c_\mathrm{AT}^2}2(\ve\|u\|^2+\ve^{-1}\|\Delta\,u\|^2)
\endaligned
$$
and minimizing with respect to $\ve$ we obtain
inequality~(\ref{Ag-Tai-vec}).
}
\end{remark}

\begin{corollary}\label{Cor-cAT}
The constant $c_\mathrm{AT}(\gamma)$
on the torus $\Omega=[0,L/\gamma]\times[0,L]$,
$\gamma\le1$,  satisfies the estimate
$c_\mathrm{AT}(\gamma)\le
 c_\mathrm{AT}/\sqrt{\gamma}\le1/\sqrt{\gamma\pi}$.
\end{corollary}
\begin{proof}
We assume for simplicity that $1/\gamma$ is integer.
Given a function $\vphi\in H^2_\mathrm{per}(\Omega)$ we extent
it by periodicity in $x_2$-direction $1/\gamma$ times
and denote the function so obtained by~$\widetilde \vphi$.
Then $\widetilde\vphi\in H^2_\mathrm{per}(\widetilde\Omega)$,
where $\widetilde\Omega=[0,L/\gamma]^2$ is a square-shaped
periodic domain so that
$\|\widetilde\vphi\|_{L_\infty(\widetilde\Omega)}\le
c_\mathrm{AT}\|\widetilde\vphi\|^{1/2}_{L_2(\widetilde\Omega)}
\|\Delta\,\widetilde\vphi\|^{1/2}_{L_2(\widetilde\Omega)}
$.
Since
$
\|\vphi\|_{L_\infty(\Omega)}=
\|\widetilde\vphi\|_{L_\infty(\widetilde\Omega)}$,
$\|\vphi\|_{L_2(\Omega)}^2=\gamma
\|\widetilde\vphi\|^2_{L_2(\widetilde\Omega)}$, and
$\|\Delta\,\vphi\|_{L_2(\Omega)}^2=\gamma
\|\Delta\,\widetilde\vphi\|^2_{L_2(\widetilde\Omega)}
$,
the corollary is proved.
\end{proof}

We now prove the remaining two inequalities used for
estimates of the number of the determining nodes.
\begin{lemma}\label{L:Jones-Titi2}{\rm (see \cite{FMRT},
\cite{Jones-Titi2}).}
Let $\Omega=[0,L_1]\times[0,L_2]$ be divided into $N$
equal squares $Q_j$ with side $l$ and let $x^j\in Q_j$
for $j=1,\dots,N$. Then for $u\in H^2_\mathrm{per}(\Omega)$ the
following inequalities hold:
\begin{align}
&\|u\|^2\le
4l^2N\,\eta^2(u)+
68l^4\|\Delta\,u\|^2,\label{Jones-Titi1}\\
&\|\nabla\,u\|^2\le2\cdot 68^{-1/2}N\eta^2(u)+68^{1/2}l^2\|\Delta\,u\|^2,
\label{Jones-Titi2}
\end{align}
where $\eta(u)=\max_{j=1,\dots,N}|u(x^j)|$.

\end{lemma}
\begin{proof}
We  consider the scalar case and prove the first inequality.
Let $u\in H^2(Q)$, where $Q=[0,l]^2$. For any
two points $\mathbf{x}=(x,y)$ and $\mathbf{x}^0=(x_0,y_0)$
in $Q$ we have
$$
u(\mathbf{x})-u(\mathbf{x}^0)=
\int_{x_0}^x u_x(\xi,y)d\xi+
\int_{y_0}^y u_y(x_0,\eta)d\eta.
$$
Hence
$$
(u(\mathbf{x})-u(\mathbf{x}^0))^2\le
2l\int_0^l u_x(\xi,y)^2d\xi+2l\int_0^l u_y(x_0,\eta)^2d\eta,
$$
which gives after integration over $Q$  with respect to $x,y$
that
$$
\|u-u(\mathbf{x}^0)\|_{L_2(Q)}^2\le
2l^2\|u_x\|_{L_2(Q)}^2+2l^3\int_0^l u_y(x_0,\eta)^2d\eta.
$$
For the second term on the right we have
$$
u_y(x_0,\eta)^2\le
u_y(x,\eta)^2+2\int_0^l|u_y(\xi,\eta||u_{yx}(\xi,\eta)|d\xi,
$$
hence, integrating with respect to $x$  and $\eta$ over $Q$ we find
$$
l\int_0^lu_y(x_0,\eta)^2d\eta\le
\|u_y\|_{L_2(Q)}^2+
2l\int_0^l\int_0^l|u_y(\xi,\eta||u_{yx}(\xi,\eta)|d\xi d\eta
\le 2\|u_y\|_{L_2(Q)}^2+l^2\|u_{xy}\|_{L_2(Q)}^2.
$$
Therefore
$$
\|u-u(\mathbf{x}^0)\|_{L_2(Q)}^2\le4l^2\|\nabla u\|_{L_2(Q)}^2+
2l^4\|u_{xy}\|_{L_2(Q)}^2.
$$
Temporarily denoting the right-side by $K$ and
using Young's inequality we have
$$
\aligned
\|u\|_{L_2(Q)}^2\le &K+2u(\mathbf{x}^0)\int_Qu(x,y)dxdy-
l^2 u(\mathbf{x}^0)^2\le\\
&K+2u(\mathbf{x}^0)\,l\|u\|_{L_2(Q)}-l^2 u(\mathbf{x}^0)^2\le
K+l^2 u(\mathbf{x}^0)^2+\frac12\|u\|_{L_2(Q)}^2.
\endaligned
$$
Hence
$$
\|u\|_{L_2(Q)}^2-2l^2 u(\mathbf{x}^0)^2\le
8l^2\|\nabla u\|_{L_2(Q)}^2+
4l^4\|u_{xy}\|_{L_2(Q)}^2.
$$
We now divide  $\Omega=[0,L_1]\times[0,L_2]$
into $N$ equal squares of side $l=(|\Omega|/N)^{1/2}$
and choose a point $\mathbf{x}^j$ in each square $Q_j$,
$j=1,\dots,N$.
Summing over $j$ we obtain
$$
\|u\|^2-2l^2\sum_{j=1}^N u(\mathbf{x}^j)^2\le
8l^2\|\nabla u\|^2+
2l^4\|\Delta\,u\|^2,
$$
where we used the periodic boundary conditions so that
 $\|u_{xy}\|^2=\int u_{xx}u_{yy}dxdy\le\frac12\|\Delta\, u\|^2$.
Next we use the interpolation inequality
$$
\|\nabla u\|^2\le\|u\|\|\Delta\,u\|\le\frac1{16l^2}\|u\|^2+
4l^2\|\Delta\,u\|^2
$$
and finally obtain
$$
\|u\|^2\le4l^2\sum_{j=1}^N u(\mathbf{x}^j)^2+
68l^4\|\Delta\,u\|^2\le
4l^2N\,\eta^2(u)+
68l^4\|\Delta\,u\|^2,
$$
which proves inequality~(\ref{Jones-Titi1})
 for the scalar case. For the vector case
we apply the above inequality for each component and
add up the results.

Finally, if $u\in H^2_\mathrm{per}(\Omega)$, $\div u=0$,
and $\rot u=\omega$,
then taking into account that $\|u\|=\|\nabla\Delta^{-1}\omega\|$
and $\|\nabla\,\omega\|=\|\Delta\,u\|$ we can write the previous
inequality in the form
\begin{equation}\label{in-vort-form}
\|\nabla\Delta^{-1}\omega\|^2\le4l^2N\eta^2(u)+
68l^4\|\nabla\,\omega\|^2=
4|\Omega|\eta^2(u)+
68|\Omega|^2N^{-2}\|\nabla\,\omega\|^2.
\end{equation}

For the proof of~(\ref{Jones-Titi2}) we have
$$
\aligned
\|\nabla\,u\|^2\le\|u\|\|\Delta\,u\|\le\
\ve\|u\|^2+(4\ve)^{-1}\|\Delta\,u\|^2\le
4Nl^2\ve\eta^2(u)+(68l^4\ve+(4\ve)^{-1})\|\Delta\,u\|^2,
\endaligned
$$
which gives~(\ref{Jones-Titi2}) by setting $\ve^{-1}=2\cdot 68^{1/2}l^2$.
\end{proof}

\section*{Acknowledgments}

 This work was supported in part by the US Civilian Research and
Development Foundation, grant no.~RUM1-2654-MO-05 ( A.A.I. and E.S.T.),
by the Russian Foundation for Fundamental Research,
grant~no.~03-01-00189 02, and by the RAS Programme
`Modern problems of theoretical
mathematics', contract no. 090703-1028 (A.A.I.).
The work of E.S.T. was supported in part by the National Science Foundation,
grant no.~DMS-0204794, the MAOF Fellowship of the  Israeli Council of Higher
Education, and
 by the USA Department of Energy under contract W-7405-ENG-36
and  the  ASCR Program in Applied Mathematical Sciences.

\bibliographystyle{amsplain}

\end{document}